# A Proof: Euler's Constant $\gamma$ is an Irrational Number

## — The answer to the front part of Hilbert's 7$^{th}$ problem

### Kaida Shi

State Key Laboratory of CAD&CG, Zhejiang University,

Hangzhou 310027, Zhejiang Province, China

Department of Mathematics, Zhejiang Ocean University,

Zhoushan 316004, Zhejiang Province, China

**Abstract**  The attributes of Euler's constant $\gamma$ have been a baffling problem to the world's mathematicians in number theory field. In 1900, when German mathematician D. Hilbert addressed the 2$^{nd}$ International Congress of Mathematicians, he suggested twenty-three previously unsolved problems to the international mathematical field. The 7$^{th}$ of these problems pertained to Euler's constant $\gamma$. After investigating this problem for many years, the author has proved that Euler's constant $\gamma$   is an irrational number.

Using the harmonic progression to proceed towards the logarithmic function, Swiss mathematician Leonard Euler suggested

$$\gamma = \lim_{n\to\infty}(1 + \frac{1}{2} + \frac{1}{3} + \Lambda + \frac{1}{n} - \log n),$$

where $\gamma$ is known as Euler's constant. For over two centuries, the attributes of Euler's constant $\gamma$ have been a focal point of the international mathematical field. In 1900, when German mathematician David Hilbert addressed the 2$^{nd}$ International Congress of Mathematicians, he suggested twenty-three mathematically baffling problems. Of these, the 7$^{th}$ contains the attribute problem of Euler's constant, namely, whether $\gamma$ is an algebraic number or a transcendental number and whether $\gamma$ is an irrational number or a rational number[1,2]? After careful investigation, I would like to express my own view of the attributes of Euler's constant $\gamma$.

## 1 The geometric implications of Euler's constant $\gamma$

If Euler's constant
$$\gamma = \lim_{n\to\infty}(1 + \frac{1}{2} + \frac{1}{3} + \Lambda + \frac{1}{n} - \log n)$$
corresponds to the geometric field, it will be:

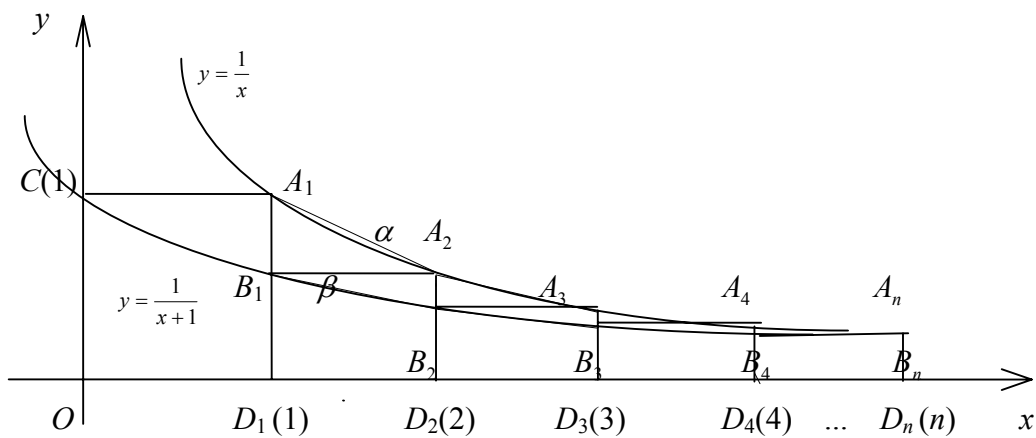

Figure

In above figure, the areas of the series of rectangles $A_1 D_1 OC, A_2 D_2 D_1 B_1, A_3 D_3 D_2 B_2, \cdots$ are respectively:
$$1, \ \frac{1}{2}, \ \frac{1}{3}, \ \Lambda \ .$$
The definite integral of the function $y = \frac{1}{x}$ on the closed interval from $D_1(1)$ to $D_n(n)$ equals the area of the curved trapezoid $A_1 A_n D_n D_1$:
$$S = \int_1^n \frac{dx}{x} = \log x \Big|_1^n = \log n \ .$$
Hence, the geometric implication of Euler's constant $\gamma$ is: when $n\to\infty$, the limit of the difference between the sum $1+\frac{1}{2}+\frac{1}{3}+\Lambda+\frac{1}{n}$ of the areas of a series of rectangles $A_1D_1OC, A_2D_2D_1B_1, A_3D_3D_2B_2, \cdots, A_nD_nD_{n-1}B_{n-1}$ and the definite integral of the

function $y=\frac{1}{x}$ on the closed interval from $D_1(1)$ to $D_n(n)$ (the area of the curved trapezoid

$A_1A_nD_nD_1$) equals $\gamma$.

## 2 The definitions of some relevant areas

（i）We take $S$ to denote the area of the **double curved trapezoid** $A_1A_nB_nB_1(n\to\infty)$. That is:

$$S' = \lim_{n\to\infty}\int_1^n (\frac{1}{x} - \frac{1}{1+x})dx$$
$$= \lim_{n\to\infty}(\log\frac{x}{1+x})_1^n$$
$$= \lim_{n\to\infty}(\log\frac{n}{1+n} - \log\frac{1}{2})$$
$$= \log 2.$$

If we divide the area $S'$ of the above double curved trapezoid as two parts $S_1$ and $S_2$, then we have: $S_1$ equals the sum of the areas of a series of curved triangles $A_1A_2B_1, A_2A_3B_2, \Lambda, A_{i-1}A_iB_{i-1}, \Lambda$. $S_2$ equals the sum of the areas of a series of curved triangles $B_1A_2B_2, B_2A_3B_3, \Lambda, B_{i-1}A_iB_i, \Lambda$.

（ii）If the areas of a series of curved triangles $A_1A_2B_1, A_2A_3B_2, \Lambda, A_{i-1}A_iB_{i-1}, \Lambda$ are

added, we have:
$$\sum_{i=2}^{\infty}\frac{1}{2}\cdot(\frac{1}{i-1}-\frac{1}{i})\cdot 1 = \sum_{i=2}^{\infty}\frac{1}{2i(i-1)} = \frac{1}{2},$$

We denote this with $S_1'$.

(iii) If the areas of a series of curved triangles $B_1A_2B_2, B_2A_3B_3, \Lambda, B_{i-1}A_iB_i, \Lambda$ are added, we have

$$\sum_{i=2}^{\infty} \frac{1}{2} \cdot (\frac{1}{i} - \frac{1}{i+1}) \cdot 1 = \sum_{i=2}^{\infty} \frac{1}{2(i+1)i} = \frac{1}{4},$$

We denote this with $S_2'$.

(iv) If we take $\alpha$ to denote the difference between $S_1'$ and $S_1$, we have

$$\alpha = S_1' - S_1 = \lim_{n\to\infty} \sum_{i=2}^{n} \left( \frac{i-\frac{1}{2}}{i(i-1)} - \log\frac{i}{i-1} \right).$$

(v) If we take $\beta$ to denote the difference between $S_2$ and $S_2'$, we have

$$\beta = S_2 - S_2' = \lim_{n\to\infty} \sum_{i=2}^{n} \left( \frac{i+\frac{1}{2}}{(i+1)i} - \log\frac{i+1}{i} \right).$$

## 3 The proof that Euler's constant $\gamma$ is irrational number

### 3.1 The proof that both $\alpha$ and $\beta$ have identical attributes

First, we will investigate the relationship between $\alpha$ and $\beta$ from a geometric angle. As long as we move a series of **little arches** (the area of the parts encircled by both broken line $B_1B_2B_3\Lambda\ B_n\Lambda$ and the curve $y = \frac{1}{1+x}$) which corresponds to the **series $\beta$** for 1 unit towards the right, they will coincide with a series of **little arches** (the area of the parts encircled by both broken line $A_2A_3A_4\Lambda\ A_{n+1}\Lambda$ and the curve $y = \frac{1}{x}$) which corresponds to the **series $\alpha$**, and the remaining part is a **little arch** encircled by both the line $A_1A_2$ and the curve $y = \frac{1}{x}$ (its area is $\frac{3}{4} - \log 2$).

If we take the front $n(n \geq 2)$ items of the following two series

$$\alpha_n = \left(\frac{3}{4} - \log 2\right) + \left(\frac{5}{12} - \log \frac{3}{2}\right) + \Lambda + \left(\frac{n-\frac{1}{2}}{n(n-1)} - \log \frac{n}{n-1}\right)$$

$$= \left(\frac{3}{4} + \frac{5}{12} + \Lambda + \frac{n-\frac{1}{2}}{n(n-1)}\right) + \left(-\log 2 - \log \frac{3}{2} - \Lambda - \log \frac{n}{n-1}\right)$$

$$= \left(\frac{3}{4} + \frac{5}{12} + \Lambda + \frac{n-\frac{1}{2}}{n(n-1)}\right) + \left(-\log 2 - \log 3 + \log 2 - \Lambda - \log n + \log(n-1)\right)$$

$$= \frac{3}{4} + \frac{5}{12} + \Lambda + \frac{n-\frac{1}{2}}{n(n-1)} - \log n;$$

and

$$\beta_n = \left(\frac{5}{12} - \log \frac{3}{2}\right) + \left(\frac{7}{24} - \log \frac{4}{3}\right) + \Lambda + \left(\frac{n+\frac{1}{2}}{(n+1)n} - \log \frac{n+1}{n}\right)$$

$$= \left(\frac{5}{12} + \frac{7}{24} + \Lambda + \frac{n+\frac{1}{2}}{(n+1)n}\right) + \left(-\log \frac{3}{2} - \log \frac{4}{3} - \Lambda - \log \frac{n+1}{n}\right)$$

$$= \left(\frac{5}{12} + \frac{7}{24} + \Lambda + \frac{n+\frac{1}{2}}{(n+1)n}\right) + \left(-\log 3 + \log 2 - \log 4 + \log 3 - \Lambda - \log(n+1) + \log n\right)$$

$$= \frac{5}{12} + \frac{7}{24} + \Lambda + \frac{n+\frac{1}{2}}{(n+1)n} + \log 2 - \log(n+1)$$

$$= \frac{5}{12} + \frac{7}{24} + \Lambda + \frac{n+\frac{1}{2}}{(n+1)n} - \log(\frac{n+1}{2}).$$

to investigate, obviously, we will discover that the series part of both series $\alpha_n$ and $\beta_n$ differ only a rational fraction, and the logarithm parts of both series $\alpha_n$ and $\beta_n$ have the identical attributes too. Therefore, $\alpha_n$ and $\beta_n$ have identical attributes.

Because

$$\alpha = \lim_{n \to \infty} \alpha_n = \lim_{n \to \infty} \left(\frac{3}{4} + \frac{5}{12} + \Lambda + \frac{n-\frac{1}{2}}{n(n-1)} - \log n\right)$$

and

$$\beta = \lim_{n\to\infty} \beta_n = \lim_{n\to\infty}\left(\frac{5}{12}+\frac{7}{24}+\Lambda+\frac{n+\frac{1}{2}}{(n+1)n}-\log(\frac{n+1}{2})\right),$$

therefore both series $\alpha$ and $\beta$ have identical attribute too. Namely, they are identically rational numbers or irrational numbers.

### 3.2 The proof that Euler's constant $\gamma$ is an irrational number

From the figures and proof results above, we have the following relationship:

$$S_1 + S_2 = \log 2, \qquad (*)$$

because

$$S_1 + \alpha = \frac{1}{2}$$

and

$$S_2 - \beta = \frac{1}{4},$$

From this, we can obtain

$$S_1 - S_2 = \frac{1}{4} - (\alpha + \beta)$$

and

$$\alpha - \beta = \frac{3}{4} - \log 2.$$

From the above relationships, we now know that the attributes of both $\alpha$ and $\beta$ will relate to the attributes of both $S_1$ and $S_2$. Now, we will focus our discussion on them.

According to the proof in Section 3.1, we know that both $\alpha$ and $\beta$ have identical attributes.

Suppose that $\alpha$ and $\beta$ are all irrational numbers, then, from $S_1 + \alpha = \frac{1}{2}$ and $S_2 - \beta = \frac{1}{4}$ we know that $S_1$ and $S_2$ are all irrational numbers. This does not contradict with the fact that $S_1 + S_2 = \log 2$ (irrational number).

Conversely, suppose that $\alpha$ and $\beta$ are all rational numbers, then, from $S_1 + \alpha = \frac{1}{2}$ and $S_2 - \beta = \frac{1}{4}$, we know that $S_1$ and $S_2$ are all rational numbers. However, because $S_1 + S_2 = \log 2$ (irrational number), these suppositions are contradictory. Hence, the possibility that $S_1$ and $S_2$ could all be rational numbers is ruled out.

Finally, because $S_1$ is an irrational number, from $S_1 = 1 - \gamma$, we know $\gamma = 1 - S_1$, Therefore, $\gamma$ is an irrational number. So, the proof is complete.

## 4 Analysis

Euler's constant $\gamma$ is not like either the constant $\pi$ or $e$. When $n \to \infty$, $\gamma$ is the limit of the difference between the harmonic series and the natural logarithm. Because Euler was alive during the 18th century, he had no advanced calculation machines to help him in his work. In the following two centuries, people could calculate at least 20,800 decimal fractions, but still did not know the attributes of Euler's constant $\gamma$. If $\gamma$ was a rational number $P/Q$, then people would have already known that $Q$ is more than 10,000 decimal fractions [3,4,5]. It is clear that even if people have advanced calculators, they do not suffice to judge the attributes of Euler's constant $\gamma$. In addition, because Euler did not consider the problem from both the analytic and geometric angles, he made the solution of his baffling problem only possible using number theory itself. This is similar to people who hope to solve 4-dimensional problems, but can only use their 3-dimensional imaginations; it brings much inconvenience and even embarrassment.